\documentclass[english,12pt,twoside]{article}
\usepackage{amssymb,amsbsy,amsmath,amsfonts,amssymb,amscd}
\usepackage{latexsym}
\usepackage{euscript}
\usepackage{exscale}
\usepackage{epsfig}
\usepackage{amsthm}
\usepackage[francais]{babel}

\newtheorem{theo}{Theorem}[section]

\newtheorem{defi}{Definition}[section]
\newtheorem{lem}{Lemma}[section]
\newtheorem{prop}{Proposition}[section]
\newtheorem{rem}{Remark}[section]

\newcommand{\Ker}{\mathrm{Ker}}

\newcommand{\St}{\mathrm St}

\newcommand{\Hom}{\mathrm{Hom}}
\newcommand{\Har}{\mathrm{Har}}

\renewcommand{\ker}{\mathrm{Ker}}

\newcommand{\Wa}{W_{\textrm{a}}}

\newcommand{\h}{\mathfrak h}
\newcommand{\al}{\alpha}

\newcommand{\ph}{\varphi}
\newcommand{\si}{\sigma}

\title{Steinberg representations and harmonic cochains for
split adjoint quasi-simple groups}
\author{Y. A\"it Amrane}

\begin{document}

\maketitle

\selectlanguage{english}
\begin{abstract}
Let $G$ be an adjoint quasi-simple group defined and split over a
non-archimedean local field $K$. We prove that the dual of the
Steinberg representation of $G$  is isomorphic to a certain space of
harmonic cochains on the Bruhat-Tits building of $G$. The Steinberg
representation is considered with coefficients in any commutative
ring.
\end{abstract}

\begin{center}
{\bf Introduction}
\end{center}

Let $K$ be a non-archimedean local field. Let $G$ be the
$K$-rational points of a reductive $K$-group of semi-simple rank
$l$. Let $T$ be a maximal $K$-split torus in $G$ and let $P$ be a
minimal parabolic $K$-subgroup of $G$ that contains $T$. There is an
abuse of language because we mean the $K$-rational points of these
algebraic subgroups of $G$. For a commutative ring $M$, the Steinberg 
representation of $G$ with coefficients in $M$ is the $M[G]$-module :
$$
\St(M)=\dfrac{C^{^\infty}(G/P,\, M)}{\sum_{Q}C^{^\infty}(G/Q,\, M)}
$$
where $Q$ runs through all the parabolic subgroups of $G$ containing
$P$.

In \cite{Borel3}, A. Borel and J.-P. Serre, computed the reduced
cohomology group $\tilde{H}^{l-1}(Y_t,M)$ of the topologized
building $Y_t$ of the parabolic subgroups of $G$ and proved that we
have an isomorphism of $M[G]$-modules :
$$
\tilde{H}^{l-1}(Y_t,M) \cong \St(M).
$$
Then they "added" this building at infinity to the Bruhat-Tits
building $X$ of $G$ to get $X$ compactified to a contractible space
$Z_t=X\amalg Y_t$. Using the cohomology exact sequence of $Z_t$ mod.
$Y_t$, they deduce an isomorphism of $M[G]$-modules :
$$
H_c^{l}(X,M) \cong \tilde{H}^{l-1}(Y_t,M).
$$
Thus,  an isomorphism of $M[G]$-modules between the compactly
supported cohomology of the Bruhat-Tits building and the Steinberg
representation of $G$ :
$$
H_c^{l}(X,M) \cong \St(M).
$$
In case $G$ is simply connected and $M$ is the complex field
$\mathbb C$, see A. Borel \cite{Borel1}, if we consider
$C^j(X,{\mathbb C})$ to be the space of $j$-dimensional cochains and
$\delta : C^{j}(X,{\mathbb C})\rightarrow C^{j-1}(X,{\mathbb C})$
the adjoint operator to the coboundary operator $d: C^{j}(X,{\mathbb
C}) \rightarrow C^{j+1}(X,{\mathbb C})$ with respect to a suitable
scalar product, we get the $l^\textrm{th}$ homology group
$H_l(X,{\mathbb C})$ of this complex as the algebraic dual of the
compactly supported cohomology group $H_c^{l}(X,{\mathbb C})$. So,
with the isomorphism above, we get a $G$-equivariant ${\mathbb
C}$-isomorphism :
$$
H_l(X,{\mathbb C})\cong \Hom_{\mathbb C}(\St({\mathbb C}),\,{\mathbb
C}).
$$
A $j$-cochain $c\in C^j(X,{\mathbb C})$ is an harmonic cochain if we
have $d(c)=\delta(c)=0$. In case of chambers $j=l$, it is clear that
we have $d(c)=0$. So if we denote by  $\Har^l({\mathbb C},{\mathbb
C})$ the space of the ${\mathbb C}$-valued harmonic cochains defined
on the chambers of $X$, we have $\Har^l({\mathbb C},{\mathbb
C})=Z_l(X,{\mathbb C})=H_l(X,{\mathbb C})$, where $Z_l(X,{\mathbb
C})=\Ker\,\delta$ is the space of the cycles at the level $l$ of the
homological complex defined by $\delta$ above. Therefore
$$
\Har^l({\mathbb C},{\mathbb C})\cong \Hom_{\mathbb C}(\St({\mathbb
C}),\,{\mathbb C}).
$$

In the present work, we consider $G$ to be a split quasi-simple
adjoint group. For any commutative ring $M$ and for any $M$-module
$L$ on which we assume $G$ acts linearly, we define $\Har^l(M,L)$ to
be the space of $L$-valued harmonic cochains on the pointed chambers
of the Bruhat-Tits building, where a pointed chamber means a chamber
with a distinguished special vertex. The notion of harmonic cochains
we use here is the same as above in case the group $G$ is also
simply connected, otherwise since we are considering pointed
chambers of the building there is an orientation property that our
cochains should also satisfy. Using a result we have proved in our
preceding paper \cite{Yacine1} that gives the Steinberg
representation of $G$ in terms of the parahoric subgroups of $G$, we
prove explicitly that we have a canonical $M[G]$-isomorphism
$$
\Har^l(M,L) \cong \Hom_M(\St(M),L).
$$
First, we give a very brief introduction to the Bruhat-Tits building
to fix our notations. Then we recall the results obtained in
\cite{Yacine1}, giving an expression of the Steinberg representation
in terms of parahoric subgroups, we will also reformulate this
result in way it becomes easier to see the link to the harmonic
cochains. Finally, we introduce the space of harmonic cochains on
the building and prove the isomorphism between this space and the
dual of the Steinberg representation of $G$.

\section{Bruhat-Tits buildings}

\paragraph{Notations}

Let $K$ be a non-archimedean local field, that is a complete field
with respect to a discrete valuation $\omega$. We assume $\omega$ to
have the value group $\omega(K^*)=\mathbb Z$.

We consider $G$ to be the group of $K$-rational points of an adjoint
quasi-simple algebraic group defined and split over $K$. Let $T$ be
a maximal split torus in $G$, $N=N_G(T)$ be the normalizer of $T$ in
$G$ and $W=N/T$ be the Weyl group of $G$ relative to $T$.

The group of characters and the group of cocharacters of $T$ are
respectively the free abelian groups
$$
X^*(T)=\Hom(T,GL_1)  \qquad\textrm{and} \qquad X_*(T)=\Hom(GL_1,T).
$$
There is a perfect duality over $\mathbb Z$
$$
\langle \cdot,\cdot  \rangle : X_*(T) \times X^*(T) \rightarrow
{\mathbb Z} \cong X^*(GL_1)
$$
with $\langle \lambda ,\chi \rangle$ given by $\chi \circ \lambda
(x)=x^{\langle \lambda ,\chi \rangle}$ for any $x\in GL_1(K)$.

Let $V=X_*(T)\otimes {\mathbb R}$  and identify its dual space $V^*$
with $X^*(T)\otimes {\mathbb R}$. Denote by $\Phi=\Phi(T,G)\subseteq
X^*(T)$ the root system of $G$ relative to $T$. By the above
duality, any root $\alpha$ induces a linear form $\alpha :V
\rightarrow {\mathbb R}$. To every root $\alpha\in \Phi$ corresponds
a coroot $\alpha^\vee \in V$, and a convolution $s_\alpha$ that acts
on $V$ by
$$
s_{\alpha}(x)=x-\langle x,\alpha \rangle \alpha^\vee.
$$
This convolution $s_\alpha$ is the orthogonal reflection with
respect to the hyperplane $H_\alpha=\ker\,\alpha$.

On the other side, we can see that the group $N$ acts on $X_*(T)$ by
conjugations. This clearly induces an action of $W$ on $V$ by linear
automorphisms. We can identify $W$ with the Weyl group $W(\Phi)$ of
the root system $\Phi$, that is the subgroup of $GL(V)$ generated by
all the reflections $s_\alpha$, $\alpha \in \Phi$.

Let $\Delta=\{1,2,\ldots , l\}$ and let $D=\{\alpha_i; i\in
\Delta\}$ be a basis of simple roots in $\Phi$. For any $i\in
\Delta$, denote $s_i=s_{\alpha_i}$. Consider $S=\{s_i;\, i\in
\Delta\}$. The pair $(W,S)$ is a finite Coxeter system.

Denote by $\Phi^\vee$ the coroot system dual to the root system
$\Phi$. Denote by $Q(\Phi^\vee)$ (resp. $P(\Phi^\vee)$) the
associated coroot lattice (resp. coweight lattice). Since we have
assumed $G$ of adjoint type we have $X_*(T)=P(\Phi^\vee)$.

\paragraph{The fundamental apartment}

Let $A_0$ be the natural affine space under $V$. Denote by
$\textrm{Aff}(A_0)$ the group of affine automorphisms of $A_0$. For
$v\in V$, denote by $\tau(v)$ the translation of $A_0$ by the vector
$v$. We have
$$
\textrm{Aff}(A_0)=V \rtimes GL(V).
$$

There is a unique homomorphism
\begin{equation}\label{nu}
\nu : T \longrightarrow X_*(T)=P(\Phi^\vee) \subseteq V
\end{equation}
such that $\langle \nu (t),\chi\rangle =-\omega(\chi(t))$ for any
$t\in T$ and any $\chi\in X^*(T)$. In our situation this
homomorphism is surjective.

An element $t\in T$ acts on $A_0$ by the translation $\tau(\nu(t))$
:
$$
tx:=\tau(\nu(t))(x)=x+\nu(t), \qquad x\in A_0,
$$
so if we put $T_0=\ker\,\nu$, this clearly induces an action of the
so-called extended affine Weyl group
$\widetilde{W}_\textrm{a}:=N/T_0$ on $A_0$. This group is an
extension of the finite group $W$ by $T/T_0$ :
$$
\widetilde{W}_\textrm{a}=\frac{N}{T_0}=\frac{T}{T_0} \rtimes W \cong
P(\Phi^\vee) \rtimes W \subseteq V \rtimes GL(V) =
\textrm{Aff}(A_0).
$$
We deduce an action of $N$ on $A_0$ by affine automorphisms that
comes from the action of $T$ by translations on $A_0$ and the linear
action of $W$ on $V$.

For any root $\alpha \in \Phi$ and any $r\in {\mathbb Z}$, let
$H_{\alpha,r}$ be the hyperplane in $A_0$ defined by
$$
H_{\alpha,r} =\{ x\in A_0;\, \langle x,\alpha \rangle - r=0 \}.
$$
Let $s_{\alpha,r}$ be the orthogonal reflection with respect to
$H_{\alpha,r}$. We have
\begin{equation}\label{reflection-translation}
s_{\alpha,r}=\tau(r\alpha^\vee) \circ s_{\alpha}.
\end{equation}
The hyperplanes $H_{\alpha,r}$ define a structure of an affine
Coxeter complex on $A_0$. Let $\Wa$ be the associated affine Weyl
group. It is a subgroup of the group $\textrm{Aff}(A_0)$ generated
by the reflections $s_{\alpha,r}$ with respect to the hyperplanes
$H_{\alpha,r}$. We have
$$
\Wa \subseteq \textrm{Aff}(A_0)=V \rtimes GL(V).
$$
In fact, $\Wa$ is the semi-direct product of $Q(\Phi^\vee)$ and $W$
(see \cite[Ch.VI,\S\,2.1,Prop. 1]{Bourbaki})
$$
\Wa = Q(\Phi^\vee) \rtimes W \subseteq P(\Phi^\vee) \rtimes W
=\widetilde{W}_\textrm{a}.
$$
The Coxeter complex $A_0$ is the fundamental apartment of the
Bruhat-Tits building.

\paragraph{The fundamental chamber}

Let $\tilde{\alpha}$ be the highest root in $\Phi$. The fundamental
chamber $C_0$ of the Bruhat-Tits building is the chamber with the
bounding walls
$$
H_{\alpha_1}=H_{\alpha_1,0}, \ldots ,H_{\alpha_l}=H_{\alpha_l,0}
\textrm{ and }H_{\tilde{\alpha},1}.
$$
It is the intersection in $A_0$ of the
open half spaces
$$
\langle x,\alpha_i \rangle > 0 \quad 1\leq i\leq l \quad\textrm{ and
}\quad \langle x,\tilde{\alpha}\rangle < 1.
$$
Denote $s_i=s_{\alpha_i}=s_{\alpha_i,0}$ for any $i$, $1\leq i\leq
l$, and $s_0=s_{\tilde{\alpha},1}$. The set $S_\textrm{a}=\{s_0,s_1,
\ldots ,s_l\}$ generates the affine Weyl group $\Wa$. The pair
$(\Wa,S_\textrm{a})$ is an affine Coxeter system and the topological
closure $\overline{C}_0$ of $C_0$ is a fundamental domain for the
action of $\Wa$ on $A_0$.

\paragraph{The Bruhat-Tits building}

The Bruhat-Tits building $X$ associated to $G$ is defined as the
quotient
$$
X=\frac{G\times A_0}{\sim}
$$
where $\sim$ is a certain equivalence relation on $G\times A_0$, see
\cite{Yacine1} or any reference on Bruhat-Tits buildings. The group
$G$ acts transitively on the chambers (the simplices of maximal
dimension) of $X$.

\section{The Steinberg representation and the Iwahori subgroup}

Let $M$ be a commutative ring on which we assume $G$ acts trivially.
For a closed subgroup $H$ of $G$, denote by $C^{^\infty}(G/H,\,M)$
(resp. $C_c^{^\infty}(G/H,\,M)$) the space of $M$-valued locally
constant functions on $G/H$ (resp. those which moreover are
compactly supported). The action of the group $G$ on the quotient
$G/H$ by left translations induces an action of $G$ on the spaces
$C^{^\infty}(G/H,\,M)$ and $C_c^{^\infty}(G/H,\,M)$.

Let $P$ be the Borel subgroup of $G$ that corresponds to the basis
$D$ of the root system $\Phi$. For any $i\in \Delta$, let
$P_i=P\coprod Ps_iP$ be the parabolic subgroup of $G$ generated by
$P$ and the reflection $s_i$. The Steinberg representation of $G$ is
the $M[G]$-module
$$
\textrm{St}(M) = \frac{C^{^\infty}(G/P,\,
M)}{\sum_{i\in\Delta}C^{^\infty}(G/P_i,M)}.
$$

Now, let $B$ be the Iwahori subgroup of $G$ corresponding to $P$.
Recall from \cite[Th. 3.4]{Yacine1} that $C^{^\infty}(G/P,\, M)$ is
generated as an $M[G]$-module by the characteristic function
$\chi_{BP}$ of the open subset $BP/P\subseteq G/P$, and then that we
have a surjective $M[G]$-homomorphism
$$
\Theta : C_c^{^\infty}(G/B,\,M) \longrightarrow C^{^\infty}(G/P,\,M)
$$
defined by $\Theta(\varphi)=\sum_{g\in G/B}\varphi(g)g.\chi_{BP}$.

For any $i\in \Delta$, let $B_i=B\coprod Bs_iB$ be the parahoric
subgroup of $G$ that corresponds to the parabolic $P_i$. Let
$\{\varpi_i;\; i\in \Delta\}$ be the fundamental coweights with
respect to the simple basis $D$ and, by the surjective homomorphism
(\ref{nu}), take $t_i\in T$ such that $\nu(t_i)=\varpi_i$. Computing
the kernel of $\Theta$, cf. [loc. cit., Th. 4.1 and Cor. 4.2], we
have :
\begin{prop}\label{STB}
We have a canonical isomorphism of $M[G]$-modules :
$$
\textrm{St}(M) \cong \frac{C_c^{^\infty}(G/B,\,
M)}{R+\sum_{i\in\Delta}C_c^{^\infty}(G/B_i,M)}
$$
where $R$ is the $M[G]$-submodule of $C_c^{^\infty}(G/B,\, M)$
generated by the functions $\chi_{Bt_iB}-\chi_{B}$, $1\leq i \leq l
$.
\end{prop}

Under the action of $G$ on the Bruhat-Tits building $X$, the Iwahori
$B$ is the pointwise stabilizer of the fundamental chamber $C_0$.
Let $B_0=B\coprod Bs_0 B$ be the parahoric subgroup of $G$ generated
by $B$ and the reflection $s_0$. The parahoric subgroups $B_i$,
$0\leq i\leq l$, are the pointwise stabilizers of the $l+1$
codimension $1$ faces of $C_0$.

We would like to reformulate the isomorphism in this proposition in
such way the connection of the Steinberg representation to harmonic
cochains on the Bruhat-Tits building looks more clear.

Denote by $l(w)$ the length of an element $w$ of the Coxeter group
$W_\textrm{a}$ with respect to the set
$S_\textrm{a}=\{s_0,s_1,\ldots ,s_l\}$ and recall that we can look
at the linear Weyl group $W$ as the subgroup of $W_\textrm{a}$
generated by the subset $S=\{s_1,\ldots,s_l\}$ of $S_\textrm{a}$.

\begin{lem}\label{bwb}
Let $g\in G$. For any $w\in W_{\rm a}$ (resp. $w\in W$), we have
$$
\chi_{BgB} -(-1)^{l(w)}\chi_{BgwB} \in
\sum^l_{i=0}C_c^{^\infty}(G/B_i,\,M) \quad \left(\textrm{resp. } \in
\sum^l_{i=1}C_c^{^\infty}(G/B_i,\,M)\right).
$$
\end{lem}
\proof Let $u_1,\ldots ,u_d \in S_\textrm{a}$ (resp. $\in S$) such
that $w=u_1 \cdots u_d$ is a reduced expression in $W_\textrm{a}$
(resp. in $W$). We have
$$
\chi_{BgB}-(-1)^d
\chi_{BgwB}=\sum_{i=1}^d(-1)^{i-1}(\chi_{Bgu_1\cdots u_{i-1}B} +
\chi_{Bgu_1\cdots u_iB}).
$$
For any $i$, if $u_i$ is the reflection $s_j$ then $
\chi_{Bgu_1\cdots u_{i-1}B} + \chi_{Bgu_1\cdots u_iB} \in
C^{^\infty}_c(G/B_j,\,M)$. \qed

Since we have assumed $G$ to be split quasi-simple, its root system
$\Phi$ is reduced and irreducible. Thus, the Dynkin diagram of the
root system $\Phi$ is one of the types described in \cite{Bourbaki},
this classification is summarized in [loc. cit., Planches I-IX].

Let $\tilde{\alpha}=\sum_{i=1}^ln_i \alpha_i$ be the highest root of
$\Phi$. From \cite[Ch.VI, \S\,2.2, Cor. of Prop. 5]{Bourbaki}, we
know that the $l+1$ vertices $v^\circ_{i}$ of the fundamental
chamber $C_0$ are $v^\circ_{0}=0$ and :
$$
v^\circ_{i}=\varpi_{i}/n_{i} \quad{\textrm for}\; 1\leq i\leq l.
$$
To each vertex $v^\circ_i$ of the fundamental chamber $C_0$ we give
the label $i$. This gives a labeling of the chamber and then of the
whole building $X$.

Denote by $J$ the subset of $\Delta=\{1,2, \ldots ,l\}$ given by
$n_i=1$. Notice that, except for a group of type $A_l$ in which all
the vertices of a chamber are special $J=\Delta$, the coroot
$\tilde{\alpha}^\vee$ dual to the highest root is equal to some
fundamental coweight $\varpi_{i_0}$, $i_0\in \Delta-J$, that induces
a special automorphism on $X$, i.e. an automorphism of $X$ that
preserves labels. So, from (\ref{reflection-translation}), we get
\begin{equation}\label{ts}
\tau(\varpi_{i_0})=\tau(\tilde{\alpha}^\vee)=s_{\tilde{\alpha},1}
s_{\tilde{\alpha}}.
\end{equation}

\begin{theo}\label{STBJ} Assume $G$ is not of type $A_l$. We have a canonical
isomorphism of $M[G]$-modules :
$$
\textrm{St}(M) \cong \frac{C_c^{^\infty}(G/B,\,
M)}{R'+\sum_{i=0}^{l} C_c^{^\infty}(G/B_i,M)}
$$
where $R'$ is the $M[G]$-submodule of $C_c^{^\infty}(G/B,\, M)$
generated by the functions $\chi_{Bt_iB}-\chi_{B}$, $i\in J$.
\end{theo}
\proof From Proposition \ref{STB}, we need to prove the equality
$$
R+\sum_{i=1}^{l} C_c^{^\infty}(G/B_i,M)=R'+\sum_{i=0}^{l}
C_c^{^\infty}(G/B_i,M).
$$

Let us prove that the left hand side is contained in the right hand
side. Let $i\in \Delta-J$. Then $t_{i}$ acts on $X$ as a special
automorphism. So, the chamber $t_{i}C_0$ is a chamber of the
apartment $A_0$ that is of the same type as $C_0$, the same type
means that any vertex $t_iv^\circ_j$ of the chamber $t_{i}C_0$ has
the same label $j$ of $v^\circ_j$. Therefore, there is $w\in
W_\textrm{a}$ such that $t_{i}C_0=wC_0$. This means that
$\chi_{Bt_{i}B}=\chi_{BwB}$ and $w$ is of even length. From Lemma
\ref{bwb}, we get
$$
\chi_{Bt_iB}-\chi_B=(-1)^{l(w)}\chi_B -\chi_B =0 \mod.
\sum_{i=0}^lC_c^{^\infty}(G/B_i, M).
$$
Therefore, $\chi_{Bt_iB}-\chi_B\in \sum_{i=0}^lC_c^{^\infty}(G/B_i,
M)$.

Now, let us prove the other inclusion. Again from Lemma \ref{bwb} we
have
$$
\chi_{B_0}=\chi_{Bs_0B}+\chi_B=\chi_{Bs_0s_{\tilde{\alpha}}s_{\tilde{\alpha}}B}+\chi_B
=(-1)^{l(s_{\tilde{\alpha}})}\chi_{Bs_0s_{\tilde{\alpha}}B}+\chi_B
\mod. \sum_{i=1}^lC_c^{^\infty}(G/B_i, M).
$$
As we have seen, (\ref{ts}), there is an $i_0\in \Delta-J$ such that
$Bt_{i_0}B=Bs_0s_{\tilde{\alpha}}B$ and $s_{\tilde{\alpha}}\in W$
being a reflection it is of odd length. Therefore,
$$
\chi_{B_0}=-\chi_{Bt_{i_0}B}+\chi_B \mod.
\sum_{i=1}^lC_c^{^\infty}(G/B_i, M),
$$
and this finishes the proof. \qed

\begin{rem}
In case $G$ is adjoint simply connected group, so of type $E_8$,
$F_4$ or $G_2$, the subset $J$ of $\Delta$ is empty, and therefore
the $M[G]$-submodule $R'$ is trivial. The theorem above gives an
isomorphism of $M[G]$-modules :
$$
\textrm{St}(M) \cong \frac{C_c^{^\infty}(G/B,\, M)}{\sum_{i=0}^{l}
C_c^{^\infty}(G/B_i,M)}.
$$
\end{rem}

\section{Steinberg representation and harmonic cochains}

Recall that the vertex $v_0^\circ$ of $C_0$ is a special vertex and
that every chamber of the building has at least one special vertex.

Let $v_i^\circ$ be a special vertex of $C_0$, this means
that $i\in J$ and that $t_i$ is a non-special automorphism of $X$.
Let $w_0$ be the longest element in $W$ and $w_{i}$ be the longest
element in the Weyl group of the root system of linear combinations
of the simple roots $\al_{j}$, $j\neq i$. Then, see \cite[Ch. VI,
\S\,2.3, Prop. 6]{Bourbaki}, we have $t_{i}w_{i}w_{0}C_0=C_0$.

Denote by $\widehat{X}^l$ the set of pointed chambers of $X$. A
pointed chamber of $X$ is a pair $(C,v)$ where $C$ is a chamber and
$v$ is a vertex of $C$ which is special. The map which to $gB$
associates the pointed chamber $g(C_0,v^\circ_0)$ gives a bijection
\begin{equation}\label{Bpointed}
G/B \xrightarrow{\sim} \widehat{X}^l.
\end{equation}
There is a natural ordering on the vertices of a pointed chamber.
Indeed, we have
$$
(C_0,v_0^\circ)=(v_0^\circ,v_1^\circ,\ldots,v_l^\circ),
$$
which corresponds to the ordering of the vertices of the extended
Dynkin diagram, and if we choose to distinguish another special
vertex in $C_0$ then the ordering on the vertices of $C_0$ will be
the one that correponds to the ordering of the vertices of the
extended Dynkin diagram we get when applying the automorphism of the
Dynkin graph that takes $0$ to the label of the new special vertex
we have chosen. We have :

\begin{lem}\label{tec}
Let $\sigma_i$ be the permutation of the set $\{0,1,\ldots ,l\}$
such that
$$
t_iw_iw_0(v^\circ_0,v^\circ_1,\ldots
,v^\circ_l)=(v^\circ_{\sigma_i(0)}, v^\circ_{\sigma_i(1)}, \ldots
,v^\circ_{\sigma_i(l)}),
$$
then
$$
{\rm sign}(\sigma_i)=(-1)^{l(w_iw_0)}.
$$
\end{lem}
\proof For any $k\in \{0,1,\ldots ,l\}$, we have
$$
v^\circ_{\sigma_i(k)}=t_iw_iw_0(v^\circ_k)=w_iw_0(v^\circ_k)+\varpi_i=w_iw_0(v^\circ_k)+v^\circ_i,
$$
thus $w_iw_0(v^\circ_k)=v^\circ_{\sigma_i(k)}-v^\circ_i$. So if we
compute the determinant of the linear automorphism $w_iw_0$ of the
vector space $V$ in the basis $(v^\circ_1, v^\circ_2,\ldots
v^\circ_{l})$, we get
$$
{\rm det}(w_iw_0)={\rm det}(v^\circ_{\sigma_i(1)}-v^\circ_i,
v^\circ_{\sigma_i(2)}-v^\circ_i, \ldots ,
v^\circ_{\sigma_i(j-1)}-v^\circ_i, -v^\circ_i,
v^\circ_{\sigma_i(j+1)}-v^\circ_i,\ldots ,
v^\circ_{\sigma_i(l)}-v^\circ_i),
$$
where $j\in \{1,2,\ldots ,l\}$ is such that $\sigma_i(j)=0$. By
subtracting the $j^\textrm{th}$ vector $-v^\circ_i$ from the other
vectors of the determinant, we get
$$
\begin{array}{ll}
{\rm det}(w_iw_0) & = - {\rm det}(v^\circ_{\sigma_i(1)},
v^\circ_{\sigma_i(2)}, \ldots, v^\circ_{\sigma_i(j-1)}, v^\circ_i,
v^\circ_{\sigma_i(j+1)},\ldots , v^\circ_{\sigma_i(l)})\\
 & = -{\rm det}(v^\circ_{\tau\sigma_i(1)},
v^\circ_{\tau\sigma_i(2)}, \ldots, v^\circ_{\tau\sigma_i(j-1)},
v^\circ_{\tau\sigma_i(j)}, v^\circ_{\tau\sigma_i(j+1)},\ldots ,
v^\circ_{\tau\sigma_i(l)})
\end{array}
$$
where $\tau=(0 \;\; i)$ is the transposition that interchanges $0$
and $i$. Therefore,
$$
{\rm det}(w_iw_0) =-{\rm sign}(\tau \sigma_i)={\rm sign}(\sigma_i),
$$
and it is clear that ${\rm det}(w_iw_0)=(-1)^{l(w_iw_0)}$.  \qed

Denote by $\hat{X}^{l-1}$ the set of all codimension one simplices
of $X$ that are ordered sets of $l$ vertices
$\eta=(v_0,\ldots,\Check{v}_i,\ldots,v_l)$ such that $v_i$ is an omitted
vertex from a pointed chamber $C=(v_0,\ldots,v_i,\ldots,v_l)\in
\hat{X}^l$. We write $\eta<C$.

Denote by $M[\hat{X}^l]$ the free $M$-module generated by the set of the pointed chambers of $X$ and let $L$ be an $M$-module on which we assume $G$ acts linearly.

\begin{defi} Let ${\mathfrak h}: M[\hat{X}^l]\rightarrow L$ be an
$M$-homomorphism. We say that $\mathfrak{h}$ is a harmonic cochain
on $X$ if it satisfies the following properties \\
{\bf (HC1)} Let $C=(v_0,v_1,\ldots ,v_l) \in \hat{X}^l$. Let
$\sigma$ be a permutation of $\{0,1,\ldots,l\}$ such that
$v_{\sigma(0)}$ is a special vertex and that
$C_\sigma=(v_{\sigma(0)},v_{\sigma(1)},\ldots ,v_{\sigma(l)})\in
\hat{X}^l$. Then
$$
{\mathfrak h}(C)=(-1)^{\textrm{sign}(\sigma)}{\mathfrak h}(C_\sigma)
$$
{\bf (HC2)} Let $\eta\in \widehat{X}^{l-1}$ be a codimension one
simplex. Let $\mathcal{B}(\eta)=\{C\in \hat{X}^l\;|\; \eta<C\}$,
then
$$
\sum_{C\in \mathcal{B}(\eta)}{\mathfrak h}(C)=0.
$$
\end{defi}

Denote by $\Har^l(M,L)$ the set of harmonic cochains.

The action of $G$ on $\Har^l(M,L)$ is induced from its natural
action on $\textrm{Hom}_M(M[\hat{X}^l],\,L)$, namely
$$
(g.\mathfrak{h})(C)=g\mathfrak{h}(g^{-1}C)
$$
for any $\mathfrak{h}\in \Har^l(M,L)$, any $g\in G$ and any $C\in
\hat{X}^l$.

\begin{rem}
In case of groups that are adjoint and simply connected, so of type
$E_8$, $F_4$ and $G_2$, there is no non-special automorphism and
therefore the first property {\bf (HC1)} of harmonic cochains is
voided.
\end{rem}

To prove the main theorem we need the following lemma 

\begin{lem}\label{bwp}
Let $g\in G$. For any $w\in W$, we have :
$$
\chi_{BgP} -(-1)^{l(w)}\chi_{BgwP} \in
\sum^l_{i=1}C^{^\infty}(G/P_i,\,M).
$$
\end{lem}
\proof The same arguments as in the proof of Lemma \ref{bwb}. \qed

\begin{theo} We have an isomorphism of $M[G]$-modules
$$
\Har^l(M,L) \cong \Hom_M(\textrm{St}(M),\, L).
$$
\end{theo}

\proof Consider the map
$$
\mathcal{H}: \Hom_M(\textrm{St}(M),\, L) \longrightarrow
\Hom_M(M[\hat{X}^l],\, L)
$$
which to $\varphi \in \Hom_M(\textrm{St}(M),\, L)$ associates
$\h_{\ph}$ defined by $\h_\ph(g(C_0,v_0^\circ))=\ph(g\chi_{BP})$ for
any $g\in G$. Let us show that $\h_\ph=\mathcal{H}(\ph)$ is a harmonic
cochain.

{\bf (HC1)} Let $v_i^\circ$ be a special vertex of $C_0$, this means that 
$v_i^\circ=t_iv_0^\circ$ with $i\in J$. Since $t_iw_iw_0$ normalizes $B$ we have 
$$
\h_\ph(C_0,v_i^\circ)=\h_\ph(t_iw_iw_0(C_0,v_0^\circ))=\ph(t_iw_iw_0\chi_{BP})=\ph(\chi_{Bt_iw_iw_0P}),
$$
and by Lemma \ref{bwp} and since $t_i\in P$, we have
$$
\ph(\chi_{Bt_iw_iw_0P})=(-1)^{l(w_iw_0)}\ph(\chi_{Bt_iP})
=(-1)^{l(w_iw_0)}\ph(\chi_{BP})=(-1)^{l(w_iw_0)}\h_\ph(C_0,v_0^\circ).
$$
Now apply Lemme \ref{tec} .

{\bf (HC2)} Let $\eta\in \hat{X}^{l-1}$. We can assume that $\eta=(v_0^\circ,v_1^\circ, \ldots
,\hat{v}_i^\circ,\ldots,v_l^\circ)$ is a face of the pointed fundamental chamber $\eta=(C_0,v_0^\circ)$. Recall from \cite{Yacine1} that  
$B_iP_i=B_iP=\coprod_{b\in B_i/B}bBP$, therefore
$$
\sum_{C\in \mathcal{B}(\eta)}\h_{\ph}(C)=\sum_{b\in
B_i/B}\h_\ph(b(C_0,v_0^\circ))=\sum_{b\in B_i/B}\ph(b\chi_{BP})=\ph(\chi_{B_iP_i})=0.
$$

Now, consider the map
$$
\Psi : \Har^l(M,L) \longrightarrow
\Hom_M(C_c^{^\infty}(G/B,\,M),\,L)
$$
which to $\h\in \Har^l(M,L)$ associates $\psi_h$ defined by
$\psi_\h(g\chi_B)=\h(g(C_0,v_0^\circ))$. Let us show that $\psi_\h$
vanishes on the $M[G]$-submodule
$R'+\sum_{i=0}^lC_c^{^\infty}(G/B_i,\,M)$ of
$C_c^{^\infty}(G/B,\,M)$. First, since $\h$ is harmonic, from {\bf
(CH2)} we deduce that for any $i$, $0\leq i\leq l$, we have
$$
\psi_\h(\chi_{B_i})=\sum_{b\in B_i/B}\psi_\h(b\chi_B)=\sum_{b\in
B_i/B}\h(b(C_0,v_0^\circ))=\sum_{C\in\mathcal{}B(\eta)}\h(C)=0,
$$
where $\eta=(v_0^\circ,\ldots ,v_{i-1}^\circ,v_{i+1}^\circ, \ldots
,v_l^\circ)$, so $\psi_\h$ vanishes on
$\sum_{i=0}^lC_c^{^\infty}(G/B_i,\,M)$. In case $t_{i}$ is a
non-special automorphism of $X$, we have
$$
\chi_{Bt_{i}B}=\chi_{Bt_{i}w_{i}w_{0}w_{0}w_{i}B}=t_{i}w_{i}w_{0}\chi_{Bw_{0}w_{i}B},
$$
therefore,
$$
\psi_{\h}(\chi_{Bt_{i}B}-\chi_{B})=\psi_{\h}((-1)^{l(w_{0}w_{i})}
t_{i}w_{i}w_{0}\chi_{B}-\chi_{B})=(-1)^{l(w_{0}w_{i})}\h(t_i w_i w_0
(C_0,v_0^\circ))-\h(C_0,v_0^\circ).
$$
Let $\si_{i}$ be the permutation of $\{0,1,\ldots ,l\}$ such that
$$
t_i w_i w_0 (v^\circ_{0},v^\circ_{1},\ldots
,v^\circ_{l})=(v^\circ_{\si_{i}(0)},v^\circ_{\si_{i}(1)}, \ldots
,v^\circ_{\si_{i}(l)}).
$$
Since $\h$ is harmonic and as so satisfy the property {\bf (HC1)},
we have
$$
\h(t_iw_iw_0(C_0,v_0^\circ))=(-1)^{\textrm{sign}(\sigma_i)}{\mathfrak
h}(C_0,v_0^\circ),
$$
Therefore,
$$
\psi_{\h}(\chi_{Bt_{i}B}-\chi_{B})=(-1)^{l(w_{0}w_{i})}(-1)^{\textrm{sign}(\sigma_i)}{\mathfrak
h}(C_0,v_0^\circ)-\h(C_0,v_0^\circ),
$$
and from Lemma \ref{tec}, we deduce that
$$
\psi_{\h}(\chi_{Bt_{i}B}-\chi_{B})=0.
$$

Finally, if we denote by $\Theta^*$ the dual homomorphism of
$\Theta$, by Theorem \ref{STBJ} we have an $M[G]$-homomorphism
$$
\Phi = {\Theta^*}^{-1}\circ \Psi : \Har^l(M,L) \longrightarrow
\Hom_M(\textrm{St}(M),\, L)
$$
which sends a harmonic cochain $\h$ to $\varphi_h$ defined by
$\varphi_\h(g\chi_{BP})=\h(g(C_0,v_0^\circ))$ for any $g\in G$. It
is easy to prove that $\Phi$ and $H$ are inverse of each other. \qed

\begin{rem}
In case $G=PGL_{l+1}(K)$, so of type $A_l$, the isomorphism in the
theorem above is established in \cite{Yacine}.
\end{rem}

\bigskip
Y. A\"{\i}t Amrane, Laboratoire Alg\`ebre et Th\'eorie des
Nombres,\\
Facult\'e de Math\'ematiques,\\
USTHB, BP 32, El-Alia, 16111 Bab-Ezzouar, Alger, Algeria. \\
e-mail : yacinait@gmail.com

\end{document}